\documentclass[11pt]{amsart}
\usepackage{amsmath,amssymb,latexsym,amsfonts,amscd}
\newtheorem{theorem}{Theorem}[section]
\newtheorem{proposition}[theorem]{Proposition}

\newtheorem{lemma}[theorem]{Lemma}

\begin{document}

\title[A property of the Frobenius map]{A property of the Frobenius map of a polynomial ring}
\author{Gennady Lyubeznik, Wenliang Zhang and Yi Zhang}
\begin{abstract} Let $R=k[x_1,\dots, x_n]$ be a ring of polynomials in a finite number of variables over  a perfect field $k$ of characteristic $p>0$  and let $F:R\to R$ be the Frobenius map of $R$, i.e. $F(r)=r^p$. We explicitly describe an $R$-module isomorphism Hom$_R(F_*(M),N)\cong {\rm Hom}_R(M,F^*(N))$ for all $R$-modules $M$ and $N$. Some recent and potential applications are discussed.\end{abstract}
\address{Dept. of Mathematics, University of Minnesota, Minneapolis,
MN 55455}
\email{gennady@math.umn.edu}
\address{Dept. of Mathematics, University of Michigan, Ann Arbor, MI 48109.}
\email{wlzhang@umich.edu}
\address{Dept. of Mathematics, University of Minnesota, Minneapolis,
MN 55455}
\email{zhang397@math.umn.edu}
\thanks{The first and third authors gratefully acknowledge NSF support through grants DMS-0202176 and DMS-0701127. The second author is supported in part by a 2009 Spring/Summer Research Fellowship from Department of Mathematics, University of Michigan.}

\maketitle

\section{Introduction}
The main result of this paper is Theorem \ref{main} which is a type of adjointness property for the Frobenius map of a polynomial ring over a perfect field. The interest in this fairly elementary result comes from its striking recent applications (see \cite[Sections 5 and 6]{WZ} and \cite{YZ}) and also from the fact that despite extensive inquiries we have not been able to find it in the published literature. 

Especially interesting is the application in \cite{YZ} where a striking new result on local cohomology modules in characteristic $p>0$ is deduced from Theorem \ref{main}. There is no doubt that the same result is true in characteristic 0, but the only currently known (to us) proof is in characteristic $p>0$, based on Theorem \ref{main}. This has motivated what promises to be a very interesting search for some new technique to extend the result in \cite{YZ} to characteristic 0.

We believe the results of this paper hold a potential for further applications and we discuss some of them in the last section.

\section{Preliminaries}
Let $R$ be a regular (but not necessarily local) UFD, let $R_s$ and $R_t$ be two copies of $R$ (the subscripts stand for {\it source} and {\it target}) and let $F:R_s\to R_t$ be a finite ring homomorphism. Let $$F_*:R_t{\rm -mod}\to R_s{\rm -mod}$$ be the restriction of scalars functor (i.e. $F_*(M)$ for every $R_t$-module $M$ is the additive group of $M$ regarded as an $R_s$-module via $F$) and let $$F^!, F^*:R_s{\rm -mod}\to R_t{\rm -mod}$$ be the functors defined by $F^!(N)={\rm Hom}_{R_s}(R_t,N)$ and $F^*(N)=R_t\otimes_{R_s}N$ for every $R_s$-module $N$. 

$F^*(N)$ and $F^!(N)$, for an $R_s$-module $N$, have a structure of $R_s$-module via the natural $R_s$-action on $R_t$, while $F_*(M)$, for an $R_t$-module $M$, retains its old $R_t$-module structure (from before the restriction of scalars). Thus $F^*(N), F^!(N)$,  and $F_*(M)$ are both $R_s$- and $R_t$-modules. 

It is well-known that $F^*$ is left-adjoint to $F_*$, i.e. there is an isomorphism of $R_t$-modules 
\begin{alignat}{2}\label{adj4}
{\rm Hom}_{R_t}(F^*(N), M)&\cong &&{\rm Hom}_{R_s}(N,F_*(M))\\
f&\mapsto &&(n\mapsto f(1\otimes n))\notag \\
(r\otimes n\mapsto rf(n))&\leftarrow &&f\notag
\end{alignat}
which is functorial in $M$ and $N$ (see, for example, \cite[II.5, p.110]{Ha}). This holds without any restrictions on $R$.

The main result of this paper is based on a different type of adjointness (see \ref{adj3} below) which is certainly not well-known (we could not find a reference). Unlike (\ref{adj4}) it is not quite canonical but depends on a choice of a certain isomorphism $\phi$ (described in (\ref{phi}) below). And it requires the conditions we imposed on $R$, namely, a regular UFD. 

Since $R_t$, being regular, is locally Cohen-Macaulay, it follows from the Auslander-Buchsbaum theorem that the projective dimension of $R_t$ as $R_s$-module is zero, i.e. $R_t$ is projective, hence locally free, as an $R_s$-module. It is a standard fact that in this case $F^!$ is right adjoint to $F_*$, \cite[Ch. 3, Exercise 6.10]{Ha} i.e. for every $R_t$-module $M$ and every $R_s$-module $N$ there is an $R_t$-module isomorphism 
\begin{alignat}{1}\label{adj1}
{\rm Hom}_{R_s}(F_*(M),N)\cong&{\rm Hom}_{R_t}(M,F^!(N))\\
(m\mapsto \mathfrak f(1))\leftarrow& f\notag
\end{alignat}
where $\mathfrak f=f(m):R_t\to N$. This isomorphism is functorial both in $M$ and in $N$.

Consider the $R_t$-module $H\stackrel{\rm def}{=}{\rm Hom}_{R_s}(R_t,R_s)$. Since $R_t$ is a locally free $R_s$-module of finite rank, it is a standard fact \cite[Ch 2, Exercise 5.1(b)]{Ha} that there is an isomorphism of functors 
\begin{align}\label{f^!}
F^!\cong H\otimes_{R_s}-,
\end{align}
i.e. for every $R_s$-module $N$ there is an $R_t$-module isomorphism 
\begin{align}\label{f^!'}
H\otimes_{R_s}N&\cong {\rm Hom}_{R_s}(R_t,N)\\
\notag h\otimes n&\mapsto (r\mapsto h(r)n)
\end{align}
and this isomorphism is functorial in $N$. Replacing $F^!$ by $H\otimes_{R_s}-$ in (\ref{adj1}) produces an $R_t$-module isomorphism 
\begin{align}\label{adj2}
{\rm Hom}_{R_s}(F_*(M),N)\cong {\rm Hom}_{R_t}(M, H\otimes_{R_s}N).
\end{align}
which is functorial in $M$ and $N$.

It follows from \cite[Kor. 5.14]{KM} that locally $H$ is the canonical module of $R_t$. In particular, the rank of $H$ as $R_t$-module is 1 and therefore $H$ is $R_t$-module isomorphic to some ideal $I$ of $R_t$. According to \cite[Kor. 6.13]{KM} the quotient $R_t/I$ is locally Gorenstein of dimension dim$R_t-1$, hence $I$ has pure height 1. Since $R_t$ is a UFD, $I$ is principal, i.e. there is an $R_t$-module isomorphism 
\begin{align}\label{phi}
\phi:R_t\to H.
\end{align} 
Hence $\phi$ induces an isomorphism of functors 
\begin{align}\label{f^*}
R_t\otimes_{R_s}-\stackrel{\phi\otimes{\rm id}}{\cong} H\otimes_{R_s}-.
\end{align}
Replacing $H\otimes_{R_s}-$ by $F^*=R_t\otimes_{R_s}-$ in (\ref{adj2}) via (\ref{f^*}) produces an $R_t$-module isomorphism 
\begin{align}\label{adj3} 
{\rm Hom}_{R_s}(F_*(M),N)\cong {\rm Hom}_{R_t}(M, F^*(N))
\end{align}
which is functorial in $M$ and $N$.

While the isomorphisms (\ref{adj1}), (\ref{f^!}), (\ref{f^!'}) and (\ref{adj2}) are canonical, the isomorphisms (\ref{phi}), (\ref{f^*}) and (\ref{adj3}) depend on a choice of $\phi$. Every $R_t$-module isomorphism $\phi':R_t\to H$ is obtained from a fixed $\phi$ by multiplication by an invertible element of $R_t$, i.e. $\phi'=c\cdot \phi$ where $c\in R_t$ is invertible. Therefore the isomorphisms (\ref{f^*}) and (\ref{adj3}) are defined up to multiplication by an invertible element of $R_t$.

Since the element $1\in R_t$ generates the $R_s$-submodule $F(R_s)$ of $R_t$ and does not belong to $\mathfrak mR_t$ for any maximal ideal $\mathfrak m$ of $R_s$, the $R_s$-module $R_t/F(R_s)$ is projective. Hence applying the functor Hom$_{R_s}(-,R_s)$ to the injective map $F:R_s\to R_t$ produces a surjection $H\to {\rm Hom}_{R_s}(R_s, R_s).$ Composing it with the standard $R_s$-module isomorphism  ${\rm Hom}_{R_s}(R_s, R_s)\stackrel{\psi\mapsto \psi(1)}{\longrightarrow}R_s$ produces an $R_s$-module surjection $H\to R_s$. Composing this latter map with the isomorphism $\phi$ from (\ref{phi}) produces an $R_s$-module surjection 
\begin{align}\label{psi}
\psi:R_t\to R_s.
\end{align} 
If $N$ is an $R_s$-module, applying $-\otimes_{R_s}N$ to $\psi$ produces an $R_s$-module surjection 
\begin{align}\label{psiN}
\psi_N:R_t\otimes_{R_s}N\to N.
\end{align}
It is not hard to check that the isomorphism (\ref{adj3}) sends $g\in {\rm Hom}_{R_t}(M,F^*(N))$ to $\psi_N\circ g\in {\rm Hom}_{R_s}(F_*(M),N).$

\section{The Main Result}
For the rest of this paper $R$ is a ring of polynomials in a finite number of variables over a perfect field $k$ of characteristic $p>0$ and $F:R_s\to R_t$ is the standard Frobenius map, i.e. $F(r)=r^p$. The main result of this paper (Theorem \ref{main}) is an explicit description of the isomorphism (\ref{adj3}) in terms of polynomial generators of $R$. The recent applications \cite{WZ, YZ} crucially depend on this explicit description. We keep the notation of the preceding section.

Let $x_1,\dots, x_n$ be some polynomial generators of $R$ over the field $k$, i.e. $R=k[x_1,\dots,x_n]$. We denote the multi-index $i_1,\dots, i_n$ by $\bar i$. Since $k$ is perfect, $R_t$ is a free $R_s$-module on the $p^n$ monomials $e_{\bar i}\stackrel{\rm def}{=}x_1^{i_1}\cdots x_n^{i_n}$ where $0\leq i_j<p$ for every $j$. If $i_j<p-1$, then $x_je_{\bar i}=e_{\bar i'}$ where $\bar i'$ is the multi-index $i_1,\dots, i_{j-1},i_j+1, i_{j+1},\dots, i_n$. If $i_j=p-1$, then $x_je_{\bar i}=x_j^pe_{\bar i'}$ where $\bar i'$ is the multi-index $i_1,\dots, i_{j-1},0,i_{j+1},\dots,i_n$.

Let $\{f_{\bar i}\in H|0\leq i_j<p\  {\rm for\ every}\  j\}$ be the dual basis of $H$, i.e. $f_{\bar i}(e_{\bar {i'}})=1$ if  $\bar i=\bar {i'}$ and $f_{\bar i}(e_{\bar {i'}})=0$ otherwise. If $i_j>0$, then $x_jf_{\bar i}=f_{\bar i'}$ where $\bar i'$ is the multi-index $i_1,\dots, i_{j-1},i_j-1,i_{j+1},\dots, i_n$. If $i_j=0$, then $x_jf_{\bar i}=x_j^pf_{\bar i'}$ where $\bar i'$ is the multi-index $i_1,\dots,i_{j-1},p-1,i_{j+1},\dots,i_n$.

Denote the multi-index $p-1,\dots, p-1$ by $\overline{p-1}$ and let $\overline{p-1}-\bar i$ be the multi-index $p-1-i_1,\dots, p-1-i_n$.

\begin{proposition}\label{Fphi}$(\rm cf.$ \cite[Remark 3.11]{BSTZ}$)$
The $R_s$-linear isomorphism $\phi:R_t\to H$ that sends $e_{\bar i}$ to $f_{\overline{p-1}-\bar i}$ is $R_t$-linear. 
\end{proposition}

\emph{Proof.} All we have to show is that $\phi(x_je_{\bar i})=x_j\phi(e_{\bar i})$ for all indices $j$ and multi-indices $\bar i$. This is straightforward from the definition of $\phi$ and the above description of the action of $x_j$ on $e_{\bar i}$ and $f_{\bar i}$.\qed

\smallskip

Clearly, $F^*(N)=R_t\otimes_{R_s}N=\oplus_{\bar i}(e_{\bar i}\otimes_{R_s}N)$, as $R_s$-modules. Thus every $R_s$-linear map $g:M\to F^*(N)$ has the form $g=\oplus_{\bar i}(e_{\bar i}\otimes_{R_s}g_{\bar i})$ where $g_{\bar i}:M\to N$ are $R_s$-linear maps (i.e. $g_{\bar i}:F_*(M)\to N$ because $M$ with its $R_s$-module structure is $F_*(M)$).

\begin{lemma} \label{F^*}
An $R_s$-linear map $g:M\to F^*(N)$ as above is $R_t$-linear if and only if $g_{\bar i}(-)=g_{\overline{p-1}}(e_{\overline{p-1}-\bar i}(-))$ for every $\bar i$ (here $e_{\overline{p-1}-\bar i}\in R_t$ acts on $(-)\in F_*(M)$ via the $R_t$-module structure on $F_*(M)$). 
\end{lemma}

\emph{Proof.} Assume $g$ is $R_t$-linear. Then $g$ commutes with multiplication by every element of $R_t$ and in particular with multiplication by $e_{\overline{p-1}-\bar i}$. That is $g(e_{\overline{p-1}-\bar i}(-))=e_{\overline{p-1}-\bar i}g(-)$. Since $e_{\overline{p-1}-\bar i}e_{\bar i}=e_{\overline{p-1}}$, the $\overline{p-1}$-component of $e_{\overline{p-1}-\bar i}g(-)$ is $e_{\overline{p-1}}\otimes_{R_s}g_{\bar i}(-)$ while the $\overline{p-1}$-component of $g(e_{\overline{p-1}-\bar i}(-))$ is $e_{\overline{p-1}}\otimes_{R_s}g_{\overline{p-1}}(e_{\overline{p-1}-\bar i}(-)).$ Since the two $\overline{p-1}$-components are equal, $g_{\bar i}(-)=g_{\overline{p-1}}(e_{\overline{p-1}-\bar i}(-))$.

Conversely, assume $g_{\bar i}(-)=g_{\overline{p-1}}(e_{\overline{p-1}-\bar i}(-))$ for every $\bar i$. To show that $g$ is $R_t$-linear all one has to show is that $g$ commutes with the action of every $x_j\in R_t$, i.e. the $\bar i$-components of $g(x_j(-))$ and $x_jg(-)$ are the same for all $\bar i$. If $i_j>0$, then the $\bar i$-component of $x_jg(-)$ is $e_{\bar i}\otimes_{R_s}g_{\bar i'}(-)$ where $\bar i'$ is the index $i_1,\dots,i_{j-1},i_j-1,i_{j+1},\dots,i_n$ while the $\bar i$-component of $g(x_j(-))$ is $e_{\bar i}\otimes_{R_s}g_{\bar i}(x_j(-))$. But the fact that $g_{\bar i}(-)=g_{\overline{p-1}}(e_{\overline{p-1}-\bar i}(-))$ for every $\bar i$ implies $g_{\bar i'}(-)=g_{\bar i}(x_j(-))$.

If $i_j=0$, then the $\bar i$-component of $x_jg(-)$ is $e_{\bar i''}{\otimes_{R_s}}(_sx_jg_{\bar i''}(-))$ where $_sx_j$ denotes the element of $R_s$ corresponding to $x_j\in R_t$ (i.e. $F(_sx_j)=x_j^p$) and $\bar i''$ is the index $i_1,\dots, i_{j-1},p-1,i_{j+1},\dots, i_n$ while the $\bar i$-component of $g(x_j(-))$ is $g_{\bar i}(x_j(-))$. But the fact that $g_{\bar i}(-)=g_{\overline{p-1}}(e_{\overline{p-1}-\bar i}(-))$ for every $\bar i$ implies $_sx_jg_{\bar i''}(-)=g_{\bar i}(x_j(-))$.\qed

\smallskip

Finally we are ready for the main result of the paper which is the following explicit description of the isomorphism (\ref{adj3}) for the Frobenius map.

\begin{theorem}\label{main}
For every $R_t$-module $M$ and every $R_s$-module $N$ there is an $R_t$-linear isomorphism 
\begin{alignat}{2}
{\rm Hom}_{R_s}(F_*(M),N)&\cong& &{\rm Hom}_{R_t}(M,F^*(N))\notag\\
g_{\overline{p-1}}(-)&\leftarrow &&(g=\oplus_{\bar i}(e_{\bar i}\otimes_{R_s}g_{\bar i}(-)))\\
g&\mapsto &&\oplus_{\bar i}(e_{\bar i}\otimes_{R_s}g(e_{\overline{p-1}-\bar i}(-))).
\end{alignat}
\end{theorem}

\emph{Proof.} As is pointed out at the end of the preceding section, the isomorphism (\ref{adj3}) sends $g\in {\rm Hom}_{R_t}(M,F^*(N))$ to $\psi_N\circ g\in {\rm Hom}_R(F_*(M),N).$ It is straightforward to check that with $\phi$ as in Proposition \ref{Fphi} the map $\psi$ of  (\ref{psi}) sends $e_{\overline{p-1}}$ to 1. This implies that $\psi_N\circ g=g_{\overline{p-1}}$ and finishes the proof that formula (11) produces the isomorphism of (\ref{adj3}) in one direction. The fact that the other direction of this isomorphism is according to formula (12) has essentially been proven in  Lemma \ref{F^*}. \qed

\section{Potential Applications}
The notion of $F$-finite modules was introduced in \cite{L}. An $F$-finite module is determined by a generating  morphism, i.e. an $R$-module homomorphism $\beta:M\to F^*(M)$ where $M$ is a finite $R$-module. For simplicity assume the $R$-module $M$ has finite length, i.e. the dimension of $M$ as a vector space over $k$, which we denote by $d$, is finite. Then the dimension of $F^*(M)$ as a vector space over $k$ equals $p^n\cdot d$. The number $p^n$ can be huge even for quite modest values of $p$ and $n$. Thus the target of  $\beta$ may be a huge-dimensional vector space even if $d, p$ and $n$ are fairly small. But the $R$-module $F_*(M)$ has dimension $d$ as a $k$-vector space, hence the map $\tilde\beta:F_*(M)\to M$ that corresponds to $\beta$ under the isomorphism of Theorem \ref{main}, is a map between two $d$-dimensional vector spaces. Huge-dimensional vector spaces do not appear! This should make the map $\tilde\beta$ easier to manage computationally than the map $\beta$. Of course the isomorphism of Theorem \ref{main} means that many properties of $\beta$ could be detected in $\tilde\beta$; for example, $\beta$ is the zero map if and only if $\tilde\beta$ is. Therein lies the potential for using Theorem \ref{main} to make computations more manageable.

The functors ${\rm Ext}^i_R(-, R)$ commute with both $F^*$ and $F_*$. More precisely, for every finitely generated $R_t$-module $M$ and every finitely generated $R_s$-module $N$  there exist functorial $R_t$-module isomorphisms $$\kappa_i:{\rm Ext}^i_{R_t}(F^*(N),R_t)\cong F^*( {\rm Ext}^i_{R_s}(N,R_s))$$$$\lambda_i:F_*({\rm Ext}^i_{R_t}(M,R_t))\cong {\rm Ext}^i_{R_s}(F_*(M),R_s).$$ Indeed, for $i=0, M=R_t$ and $N=R_s$ a straightforward composition of the $R_t$-module isomorphism $R_t\otimes_{R_s}R_s\stackrel{r'\otimes r\mapsto r'r^p}{\longrightarrow}R_t$ and the standard $R$-module isomorphisms ${\rm Hom}_{R}(R,R)\cong R$ for $R=R_t, R_s$ produce $\kappa_0$ while an additional $R_t$-module isomorphism ${\rm Hom}_{R_s}(R_t,R_s)=H\cong R_t$ produces $\lambda_0$. This implies that if $-$ stands for a complex of finite free $R$-modules, then ${\rm Hom}_R(-,R)$ commutes with $F^*$ and $F_*$. Taking now finite free resolutions of $M$ and $N$ in the categories of $R_t$- and $R_s$-modules respectively and considering that $F^*$ and $F_*$, being exact, commute with the operation of taking the (co)homology of complexes, we get $\kappa_i$ and $\lambda_i$ for every $i$.

It is straightforward to check that there is a commutative diagram
$$
\begin{CD}
{\rm Hom}_{R_t}(F^*(N), M)@>>> {\rm Hom}_{R_s}(N,F_*(M))\\
@VVV @VVV\\
{\rm Hom}_{R_t}({\rm Ext}^i_{R_t}(M,R_t),{\rm Ext}^i_{R_t}(F^*(N),R_t))@.{\rm Hom}_{R_s}({\rm Ext}^i_{R_s}(F_*(M), R_s),{\rm Ext}^i_{R_s}(N,R_s))\\
@VVV @VVV\\
{\rm Hom}_{R_t}({\rm Ext}^i_{R_t}(M,R_t),F^*({\rm Ext}^i_{R_s}(N,R_s)))@>>>{\rm Hom}_{R_s}(F_*({\rm Ext}^i_{R_t}(M, R_t)),{\rm Ext}^i_{R_s}(N,R_s))
\end{CD}
$$
where the top horizontal map is the isomorphism (\ref{adj4}), the bottom horizontal map is the isomorphism of Theorem \ref{main}, the bottom vertical maps are induced by $\kappa_i$ and $\lambda_i$ and, finally, the top vertical maps are defined by sending every $f\in {\rm Hom}(L,L')$ to the map ${\rm Ext}^i(L', R)\to {\rm Ext}^i(L,R)$ functorially induced by $f$. 

In other words, if a pair of maps $F^*(N)\to M$ and $N\to F_*(M)$ correspond to each other under (\ref{adj4}), then the induced maps ${\rm Ext}^i_{R_t}(M,R_t)\to F^*({\rm Ext}^i_{R_s}(N,R_s))$ and $F_*({\rm Ext}^i_{R_t}(M, R_t)) \to {\rm Ext}^i_{R_s}(N,R_s)$ correspond to each other under isomorphism of Theorem \ref{main}.

Of course there is a similar diagram with the isomorphism of Theorem \ref{main} in the top row; we are not going to use it in the rest of the paper. 

\smallskip

An important example of $F$-finite modules are local cohomology modules $H^i_I(R)$ of $R$ with support in an ideal $I\subset R$. A generating morphism of $H^i_I(R)$ is  the composition $$f:{\rm Ext}^i_R(R/I, R)\to {\rm Ext}^i_R(F^*(R/I), R)\stackrel{\kappa_i}{\cong}F^*({\rm Ext}^i_R(R/I, R))$$
where the first map is induced by the isomorphism $F^*(R/I)\stackrel{r'\otimes\bar r\mapsto r'\bar r^p}{\cong}R/I^{[p]}$ followed by the natural surjecton $R/I^{[p]}\to R/I$. The following proposition holds the potential for simplifying computations involving local cohomology modules.

\begin{proposition}
Let $I\subset R$ be an ideal and let the composition $$f:{\rm Ext}^i_R(R/I, R)\to {\rm Ext}^i_R(F^*(R/I), R)\stackrel{\kappa_i}{\cong}F^*({\rm Ext}^i_R(R/I, R))$$ be as above.The map that corresponds to $f$ under the isomorphism of Theorem \ref{main} is the composition $$g:F_*({\rm Ext}^i_R(R/I, R))\stackrel{\lambda_i}{\cong}{\rm Ext}^i_R(F_*(R/I), R)\to {\rm Ext}^i_R(R/I, R)$$ where the second map in the composition is nothing but the map induced on ${\rm Ext}^i_R(-,R)$ by the natural Frobenius map $R/I\stackrel{r\mapsto r^p}{\cong} F_*(R/I)$. 
\end{proposition}

\emph{Proof.} This is immediate from the above commutative diagram considering that  the maps $F^*(R/I)\stackrel{r'\otimes r\mapsto r'r^p}{\longrightarrow} R/I$ and $R/I\stackrel{r\to r^p}{\rightarrow}F_*(R/I)$ correspond to each other under the isomorphism (\ref{adj4}).\qed

\smallskip

In addition to potential use in computation, the material of this section has already been used in a proof of a theoretical result \cite[Section 5]{WZ}.

\end{document}